# The symmetry of the Kepler problem, the inverse Ligon-Schaaf mapping and the Birkhoff conjecture

**Thomas Sumner Ligon\***

Faculty of Physics and Center for NanoScience (CeNS), Ludwig-Maximilians-Universität, München, Germany

\* Thomas.Ligon@physik.uni-muenchen.de

## Abstract

The Ligon-Schaaf regularization (*LS* mapping) was introduced in 1976 and has been used several times. However, we are not aware of any direct usage of the inverse mapping, perhaps since it appears at first sight to be quite complex, involves the use of a transcendental equation (referred to as the generalized Kepler equation) that cannot be solved in closed form, and lacks smoothness near the collision point. Here, we provide some insight into the significance of this equation, along with a very simple derivation and confirmation of the inverse *LS* mapping. Then we use numerical methods to investigate three applications: 1) solutions of the Kepler function, 2) calculation of orbits including time-of-flight data based on the Delaunay Hamiltonian, and 3) numerical evidence for the Birkhoff conjecture for the circular restricted 3-body problem.

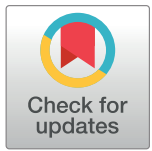

**Data Availability Statement:** All relevant data are synthetic data created by software. Both software and data are included in the Supporting Information files.

**Funding:** The author received no specific funding for this work.

**Competing interests:** The author has declared that no competing interests exist.

## Introduction

The Ligon-Schaaf mapping (*LS* mapping) was introduced in 1976 [1] based in part on a Diplom thesis (similar to master's thesis) which has recently been translated into English [2–4]. Since then, the paper has been cited numerous times [5–27]. In particular, a few papers have provided very significant insight into the properties of the mapping and the symplectic manifolds involved in the regularization that it achieves [8, 17].

The *LS* mapping is a modified stereographic projection, converting phase space in Cartesian coordinates to phase space over a 3-sphere, and maps the collision orbits to the north pole of the sphere. This makes it possible to regularize the Kepler problem by adding the poles to the original space. For the analysis of symmetry, it maps both angular momentum and the Runge-Lenz vector to angular momentum on the sphere and maps the Kepler Hamiltonian to the Delaunay Hamiltonian. The result is a symmetry-preserving diffeomorphism between the two phase spaces.

The thesis which preceded the paper had the goal of calculating the classical Kepler problem as an example of how symplectic differential geometry could be used in theoretical physics. After introducing the concepts and discussing some traditional topics, such as the conservation of angular momentum, we examined the conservation of the Runge-Lenz vector and asked ourselves if the corresponding symmetry was global. In chapter VIII, we searched for the integral curves of the Runge-Lenz vector by calculating the Poisson brackets with the Cartesian coordinates $p$, $q$. However, these equations were quite complex, and we didn't find a full solution. As a kind of coordinate transformation, in chapter IX, we applied the Moser mapping [28], which is based on the stereographic projection that had been used for the quantum-mechanical analysis of the hydrogen atom [29]. After applying this mapping, the flow of the Runge-Lenz vector could be integrated directly, demonstrating that this is in fact a global symmetry. However, the solutions involved two different angles, and the second one required solving a transcendental equation that cannot be solved in closed form. Then we discovered that we could modify the Moser mapping in such a way that this transcendental equation disappeared from the solutions, and both angular momentum and the Runge-Lenz vector generated simple rotations on the sphere $S^3$. In other words, we had defined the *LS* mapping and showed that it converts the symmetry of the Kepler problem to the canonical action of SO(4).

In the paper of 1976, the *LS* mapping was presented in detail, the symmetry properties of the symplectic manifolds were discussed, and the inverse mapping was presented. Here, the transcendental equation that we removed from the action of the Runge-Lenz vector reappears as part of the inverse mapping, and it was shown to be a kind of "generalized Kepler equation".

The inverse *LS* mapping lacks smoothness in the sense that, given a smooth function $f(q,p)$ or even just $f(q)$, the function $f \circ \Phi_{LS}^{-1} = f \circ F^{-1}$ does not necessarily extend to a smooth function on the full regularized set, which is the cotangent bundle of a sphere. This has been a problem in several applications, as discussed for example in the context of the averaging method in [8] (last paragraph before section 2). We also expect that this will cause problems in the context of Kolmogorov-Arnold-Moser (KAM) theory. We observed this lack of smoothness near the collision point in Observation 10 (see S1 Document, Discussion of Observation 10 and S6 Fig. Sun-Jupiter energy and curvature 2).

In this paper, we apply the factorization of the *LS* mapping, as formulated by Cushman and Bates [9]. This way, the mapping is split into an algebraic part and a trigonometric part, which is basically a rescaled rotation. Both parts are very easy to invert, which we do in the main text of the paper. In particular, the inverse of the rescaled rotation is another rotation in the opposite direction, which shows that the angle used in the inverse mapping is exactly the same angle as in the forward mapping. When we take account of this fact, all calculations involving the inverse mapping are greatly simplified. Then we calculate some identities and observe that the angle is transformed this way: $\varphi = \sqrt{-2H}\boldsymbol{q} \cdot \boldsymbol{p} = \xi_0 sin(\varphi) - \frac{\eta_0}{\eta} cos(\varphi)$. In other words, the transcendental equation used by the inverse mapping is the only way we have of finding the angle $\varphi$ when we move in the reverse direction.

Based on this foundation, we create a numeric solution for all mappings and use it for the following applications: 1) investigation of the "Kepler function" generated by the "generalized Kepler equation": $\varphi = x \cdot sin(\varphi) - y \cdot cos(\varphi)$, 2) calculation of time-of-flight for Kepler orbits using the Delaunay Hamiltonian and the inverse *LS* mapping, and 3) numerical evidence for the Birkhoff conjecture for the circular restricted 3-body problem. It has often been observed that the Kepler laws can be deduced by geometrical means, without solving differential equations, but finding the time of flight requires solving the Kepler equation. In our second application, this time dependency is found in closed form on $T^*S^3$, and the Kepler equation is solved "implicitly" as part of the inverse *LS* mapping. The third application provides evidence for the Birkhoff conjecture within the framework of holomorphic curve theory, but all calculations in this paper are "elementary calculations" based solely on traditional algebra and trigonometry. This investigation begins with an alternative, analytical proof for the Birkhoff conjecture based on the rotating Kepler Hamiltonian[13] and then covers numerical calculations based on the circular restricted 3-body Hamiltonian.

## Negative energy

### Notation

We begin with an overview and a comparison with notation used by Cushman and Bates[9]. Fig 1 provides an overview of the mapping in our notation.

Here, we use the notation from our original paper, and provide a detailed comparison with the Cushman and Bates notation in the supporting information.

### Forward LS mapping $F = F_2 \circ F_1$

The forward *LS* mapping is defined for $H < 0$, corresponding to the elliptical orbits of the Kepler Hamiltonian, and $q \neq 0$, meaning that collision orbits are excluded.

**Definition 1.** *For negative energy and $q \neq 0$, the forward LS mapping is defined by*

i. *Mapping **F***

$$H = \frac{p^2}{2} - \frac{1}{q} \tag{1}$$

$$\varphi = \sqrt{-2H} \boldsymbol{q} \cdot \boldsymbol{p} \tag{2}$$

$$\xi_0 = (p^2 q - 1)cos(\varphi) + \sqrt{-2H} \boldsymbol{q} \cdot \boldsymbol{p}\ sin(\varphi) \tag{3}$$

$$\boldsymbol{\xi} = \sqrt{-2H} q \boldsymbol{p} cos(\varphi) + \left[\frac{\boldsymbol{q}}{q} - (\boldsymbol{q} \cdot \boldsymbol{p})\boldsymbol{p}\right] sin(\varphi) \tag{4}$$

$$\eta_0 = -\boldsymbol{q} \cdot \boldsymbol{p}\ cos(\varphi) + \frac{1}{\sqrt{-2H}}(p^2 q - 1)sin(\varphi) \tag{5}$$

$$\boldsymbol{\eta} = -\frac{1}{\sqrt{-2H}}\left[\frac{\boldsymbol{q}}{q} - (\boldsymbol{q} \cdot \boldsymbol{p})\boldsymbol{p}\right]cos(\varphi) + q \boldsymbol{p} sin(\varphi) \tag{6}$$

ii. *Mapping **F_1***

$$r_0 = p^2 q - 1 \tag{7}$$

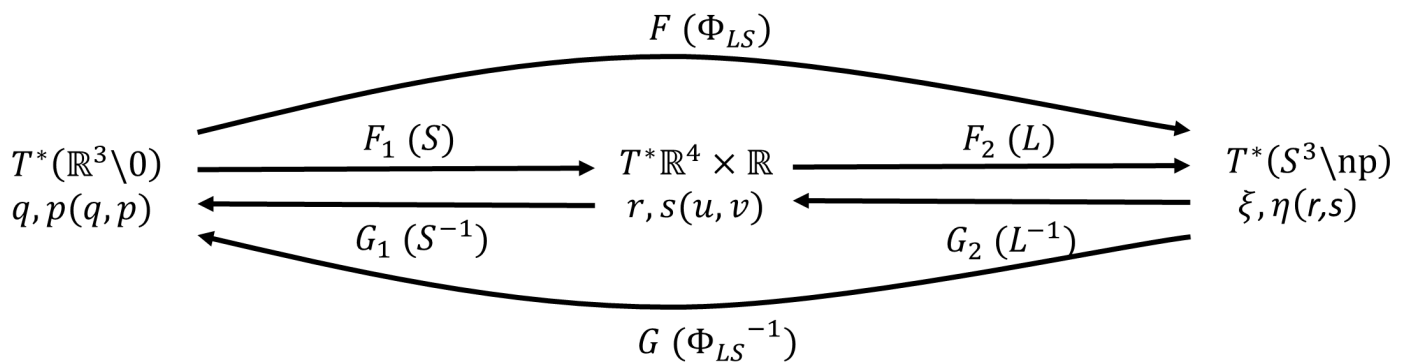

**Fig 1. LS mapping.** The *LS* mapping (forward mapping $F$ and inverse mapping $G$). The factorization consists of an algebraic part ($F_1$) and a trigonometric part ($F_2$). The notation of Cushman and Bates is shown in parentheses.

$$\boldsymbol{r} = \sqrt{-2H}q\boldsymbol{p} \tag{8}$$

$$s_0 = -\sqrt{-2H}\boldsymbol{q} \cdot \boldsymbol{p} \tag{9}$$

$$\boldsymbol{s} = -\left[\frac{\boldsymbol{q}}{q} - (\boldsymbol{q} \cdot \boldsymbol{p})\boldsymbol{p}\right] \tag{10}$$

iii. *Mapping $F_2$*

$$\varphi = -s_0 \tag{11}$$

$$F_2 = \begin{pmatrix} cos(\varphi) & -sin(\varphi) \\ \dfrac{1}{\sqrt{-2H}}sin(\varphi) & \dfrac{1}{\sqrt{-2H}}cos(\varphi) \end{pmatrix} \tag{12}$$

$$\xi_0 = r_0 cos(\varphi) - s_0 sin(\varphi) \tag{13}$$

$$\boldsymbol{\xi} = \boldsymbol{r}cos(\varphi) - \boldsymbol{s}\,sin(\varphi) \tag{14}$$

$$\eta_0 = \frac{1}{\sqrt{-2H}}s_0 cos(\varphi) + \frac{1}{\sqrt{-2H}}r_0 sin(\varphi) \tag{15}$$

$$\boldsymbol{\eta} = \frac{1}{\sqrt{-2H}}\boldsymbol{s}\,cos(\varphi) + \frac{1}{\sqrt{-2H}}\boldsymbol{r}\,sin(\varphi) \tag{16}$$

The mappings $\Phi_{LS} = L \circ S$ in the Cushman and Bates notation, corresponding to $F = F_2 \circ F_1$, are provided in the supporting information.

***Proposition 1.*** $F = F_2 \circ F_1$.

**Proof.** This is a simple substitution. ∎

In the factorization presented by Cushman and Bates, $\boldsymbol{F_1}$ is the algebraic part and $\boldsymbol{F_2}$ is the trigonometric part, which can be understood as a rescaled rotation.

At this point, the mapping $F_1$ can be inverted directly (giving $G_1$), via elementary calculations (see below), and the matrix $F_2$ can be inverted directly (giving $G_2$). Then, they can be combined to give $G$. The definition of $\varphi$ follows from identity (36).

***Definition 2.*** *For negative energy, the inverse LS mapping is defined by*

i. *Mapping $G$*

$$\varphi = \xi_0 sin(\varphi) - \frac{\eta_0}{\eta}cos(\varphi) \tag{17}$$

$$\boldsymbol{q} = \eta^2 \left[ -\boldsymbol{\xi} \left( \frac{\eta_0}{\eta} - sin(\varphi) \right) + \frac{\boldsymbol{\eta}}{\eta} \left( \xi_0 - cos(\varphi) \right) \right] \tag{18}$$

$$\boldsymbol{p} = \frac{\boldsymbol{\xi} cos(\varphi) + \frac{\boldsymbol{\eta}}{\eta} sin(\varphi)}{\eta \left( 1 - \xi_0 cos(\varphi) - \frac{\eta_0}{\eta} sin(\varphi) \right)} \tag{19}$$

ii. *Mapping $G_1$*

$$\boldsymbol{q} = \frac{1}{2H} \left[ \boldsymbol{s}(1 - r_0) + \boldsymbol{r} s_0 \right] \tag{20}$$

$$\boldsymbol{p} = \frac{\sqrt{-2H} \boldsymbol{r}}{(1 - r_0)} \tag{21}$$

iii. *Mapping $G_2$*

$$\varphi = \xi_0 sin(\varphi) - \frac{\eta_0}{\eta} cos(\varphi) \tag{22}$$

$$G_2 = \begin{pmatrix} cos(\varphi) & \sqrt{-2H} sin(\varphi) \\ -sin(\varphi) & \sqrt{-2H} cos(\varphi) \end{pmatrix} \tag{23}$$

$$r_0 = \xi_0 cos(\varphi) + \frac{\eta_0}{\eta} sin(\varphi) \tag{24}$$

$$\boldsymbol{r} = \boldsymbol{\xi} cos(\varphi) + \frac{\boldsymbol{\eta}}{\eta} sin(\varphi) \tag{25}$$

$$s_0 = -\xi_0 sin(\varphi) + \frac{\eta_0}{\eta} cos(\varphi) \tag{26}$$

$$\boldsymbol{s} = -\boldsymbol{\xi} sin(\varphi) + \frac{\boldsymbol{\eta}}{\eta} cos(\varphi) \tag{27}$$

The value of $\varphi$ in (17) comes from the identity (36).

The denominator of (19), $1 - \xi_0 cos(\varphi) - \frac{\eta_0}{\eta} sin(\varphi) \neq 0$, and the denominator of (21), $1 - r_0 \neq 0$, because, following the identity (37), $-2Hq = 1 - r_0 = 1 - \xi_0 cos(\varphi) - \frac{\eta_0}{\eta} sin(\varphi)$, and the mapping is only defined for $H < 0$, corresponding to the elliptical orbits of the Kepler Hamiltonian, and $q \neq 0$, meaning that collision orbits are excluded.

Proposition 2.

i. $G = G_1 \circ G_2$.

ii. $G_1 = F_1^{-1}$

iii. $G_2 = F_2^{-1}$

iv. $G = F^{-1}$.

**Proof.** (i) We begin with the definition of $\boldsymbol{q}$ from $G_1$,

$$\boldsymbol{q} = \frac{1}{2H}[s(1 - r_0) + \boldsymbol{r}s_0]$$

substitute the values for $\boldsymbol{r}$ and $\boldsymbol{s}$ from $G_2$,

$$\boldsymbol{q} = \frac{1}{2H}\left[\left(-\xi sin(\varphi) + \frac{\boldsymbol{\eta}}{\eta}cos(\varphi)\right)\left(1 - \xi_0 cos(\varphi) - \frac{\eta_0}{\eta}sin(\varphi)\right) + \left(\xi cos(\varphi) + \frac{\boldsymbol{\eta}}{\eta}sin(\varphi)\right)\left(-\xi_0 sin(\varphi) + \frac{\eta_0}{\eta}cos(\varphi)\right)\right]$$

expand the expression,

$$\boldsymbol{q} = \frac{1}{2H}\left[-\xi sin(\varphi) + \frac{\boldsymbol{\eta}}{\eta}cos(\varphi) + \xi\xi_0 sin(\varphi)cos(\varphi) - \xi_0\frac{\boldsymbol{\eta}}{\eta}cos^2(\varphi) + \xi\frac{\eta_0}{\eta}sin^2(\varphi)\right.$$

$$\left. -\frac{\eta_0}{\eta^2}sin(\varphi)cos(\varphi) - \xi\xi_0 sin(\varphi)cos(\varphi) - \xi_0\frac{\eta}{\eta}sin^2(\varphi) + \xi\frac{\eta_0}{\eta}cos^2(\varphi) + \frac{\eta_0}{\eta^2}sin(\varphi)cos(\varphi)\right]$$

collect similar terms,

$$\boldsymbol{q} = \frac{1}{2H}\left[-\xi sin(\varphi) + \frac{\boldsymbol{\eta}}{\eta}cos(\varphi) - \xi_0\frac{\boldsymbol{\eta}}{\eta} + \xi\frac{\eta_0}{\eta}\right]$$

and substitute the value of $H$ from (35).

$$\boldsymbol{q} = \eta^2\left[-\xi\left(\frac{\eta_0}{\eta} - sin(\varphi)\right) + \frac{\boldsymbol{\eta}}{\eta}(\xi_0 - cos(\varphi))\right]$$

Now we take the definition of $\boldsymbol{p}$ from $G_1$.

$$\boldsymbol{p} = \frac{\sqrt{-2H}\boldsymbol{r}}{(1 - r_0)}$$

and substitute the values or $\boldsymbol{r}$ and $\boldsymbol{s}$ from $G_2$.

$$\boldsymbol{p} = \frac{\left(\xi cos(\varphi) + \frac{\eta}{\eta}sin(\varphi)\right)}{\eta\left(1 - \xi_0 cos(\varphi) - \frac{\eta_0}{\eta}sin(\varphi)\right)}$$

(ii) We begin with $F_1$ (Eqs (7) through (10)).
Now we can solve (7) for $q$:

$$r_0 - 1 = p^2 q - 2 = 2q\left(\frac{p^2}{2} - \frac{1}{q}\right) = 2Hq.$$

$$q = \frac{1}{-2H}(1 - r_0) \tag{28}$$

Then we can solve (8) for $\boldsymbol{p}$:

$$\boldsymbol{p} = \frac{\boldsymbol{r}}{\sqrt{-2Hq}} = \sqrt{-2H}\frac{\boldsymbol{r}}{(1 - r_0)}.$$

This gives us Eq (21). With this, we can solve (10) for $\boldsymbol{q}$:

$$\boldsymbol{q} = -q\boldsymbol{s} + (\boldsymbol{q} \cdot \boldsymbol{p})q\boldsymbol{p} \tag{29}$$

substitute $q$, $\boldsymbol{q} \cdot \boldsymbol{p}$ and $\boldsymbol{p}$ from (28), (9), and (21),

$$\boldsymbol{q} = \frac{1}{2H}(1 - r_0)\boldsymbol{s} + \frac{s_0}{-\sqrt{-2H}}\left(\frac{1}{-2H}\right)(1 - r_0)\sqrt{-2H}\frac{\boldsymbol{r}}{(1 - r_0)}$$

and simplify.

This gives us Eq (20). The result is now mapping $G_1$.

(iii) The mapping $G_2$, the inverse of $F_2$, can be obtained by reversing the direction of rotation and inverting the scaling factor. In addition, we also use a formula for the angle $\varphi$ (36) and for expressing the Hamiltonian in terms of the Delaunay Hamiltonian $-\frac{1}{2\eta^2}$ in (35).

(iv) This is a direct result of (i)-(iii). ∎

The corresponding derivation of $\Phi_{LS}^{-1} = S^{-1} \circ L^{-1}$ by substitution and the mappings $\Phi_{LS}^{-1} = S^{-1} \circ L^{-1}$, $L^{-1}$ and $S^{-1}$ in the Cushman and Bates notation are provided in the supporting information.

## Identities

**Proposition 3.** *The following identities hold*:

$$r^2 = \sum_{i=0}^{3} r_i^2 = 1 \tag{30}$$

$$\sum_{i=0}^{3} r_i s_i = 0 \tag{31}$$

$$s^2 = \sum_{i=0}^{3} s_i^2 = 1 \tag{32}$$

$$\xi^2 = \sum_{i=0}^{3} \xi_i^2 = 1 \tag{33}$$

$$\sum_{i=0}^{3} \xi_i \eta_i = 0 \tag{34}$$

$$\eta^2 = \sum_{i=0}^{3} \eta_i^2 = -\frac{1}{2H} \tag{35}$$

$$\varphi = \sqrt{-2H}\boldsymbol{q} \cdot \boldsymbol{p} = -s_0 = \xi_0 sin(\varphi) - \frac{\eta_0}{\eta}cos(\varphi) \tag{36}$$

$$-2Hq = -p^2 q + 2 = 1 - r_0 = 1 - \xi_0 cos(\varphi) - \frac{\eta_0}{\eta}sin(\varphi) \tag{37}$$

$$L_i = q_j p_k - q_k p_j = \frac{1}{\sqrt{-2H}}\left(r_j s_k - r_k s_j\right) = \xi_j \eta_k - \xi_k \eta_j \tag{38}$$

$$M_i = \frac{1}{\sqrt{-2H}}\left(\frac{q_i}{q} + p_i(\boldsymbol{q} \cdot \boldsymbol{p}) - q_i p^2\right) = \frac{1}{\sqrt{-2H}}(r_0 s_i - r_i s_0) = \xi_0 \eta_i - \xi_i \eta_0 \tag{39}$$

$$H = \frac{p^2}{2} - \frac{1}{q} = -\frac{1}{2\eta^2} \tag{40}$$

$$\frac{p^2}{2} - \frac{1}{q} + q_2 p_1 - q_1 p_2 = H + \frac{1}{\sqrt{-2H}} (r_2 s_1 - r_1 s_2) = -\frac{1}{2\eta^2} + \xi_2 \eta_1 - \xi_1 \eta_2 \qquad (41)$$

**Proof.** These are all direct results the definitions and can be calculated in either direction, i.e. solely based on the forward mapping or the inverse mapping (details in supporting information). ∎

The identity (36) for the angle $\varphi$ is very important, because it tells us how to calculate $\varphi$ in the inverse mapping. We refer to the resulting equation as the "generalized Kepler equation".

In identities (38) and (39), the angular momentum and the rescaled Runge-Lenz vector are mapped to angular momentum. Identity (40) maps the Kepler Hamiltonian to the Delaunay Hamiltonian. Identity (41) shows that the Hamiltonian of the rotating Kepler problem also keeps its form.

The Hamiltonian of the circular restricted 3-body problem does not yield a simple form; the terms involving $sin(\varphi)$ and $cos(\varphi)$ do not cancel out. We treat this case numerically in the applications.

## Positive energy

### Mapping $F = F_2 \circ F_1$ (LS forward)

This is the forward *LS* mapping, and is defined for $H > 0$, corresponding to the hyperbolic orbits of the Kepler Hamiltonian, and $q \neq 0$, meaning that collision orbits are excluded.

**Definition 3.** *For positive energy and $q \neq 0$, the forward LS mapping is defined by*

i. *Mapping $F$ (the first equation is the same as Eq (1))*

$$H = \frac{p^2}{2} - \frac{1}{q}$$

$$\varphi = \sqrt{2H} \boldsymbol{q} \cdot \boldsymbol{p} \qquad (42)$$

$$\xi_0 = (p^2 q - 1) cosh(\varphi) - \sqrt{2H} \boldsymbol{q} \cdot \boldsymbol{p}\, sinh(\varphi) \qquad (43)$$

$$\boldsymbol{\xi} = \sqrt{2H} q \boldsymbol{p}\, cosh(\varphi) + \left[ \frac{\boldsymbol{q}}{q} - (\boldsymbol{q} \cdot \boldsymbol{p}) \boldsymbol{p} \right] sinh(\varphi) \qquad (44)$$

$$\eta_0 = \boldsymbol{q} \cdot \boldsymbol{p}\, cosh(\varphi) - \frac{1}{\sqrt{2H}} (p^2 q - 1) sinh(\varphi) \qquad (45)$$

$$\boldsymbol{\eta} = -\frac{1}{\sqrt{2H}} \left[ \frac{\boldsymbol{q}}{q} - (\boldsymbol{q} \cdot \boldsymbol{p}) \boldsymbol{p} \right] cosh(\varphi) - q \boldsymbol{p}\, sinh(\varphi) \qquad (46)$$

ii. *Mapping $F_1$*

$$r_0 = p^2 q - 1 \qquad (47)$$

$$\boldsymbol{r} = \sqrt{2H} q \boldsymbol{p} \qquad (48)$$

$$s_0 = -\sqrt{2H}\boldsymbol{q} \cdot \boldsymbol{p} \tag{49}$$

$$\boldsymbol{s} = \left[\frac{\boldsymbol{q}}{q} - (\boldsymbol{q} \cdot \boldsymbol{p})\boldsymbol{p}\right] \tag{50}$$

iii. *Mapping $F_2$*

$$\varphi = -s_0 \tag{51}$$

$$F_2 = \begin{pmatrix} cosh(\varphi) & sinh(\varphi) \\ -\dfrac{1}{\sqrt{2H}}sinh(\varphi) & -\dfrac{1}{\sqrt{2H}}cosh(\varphi) \end{pmatrix} \tag{52}$$

$$\xi_0 = r_0 cosh(\varphi) + s_0 sinh(\varphi) \tag{53}$$

$$\boldsymbol{\xi} = \boldsymbol{r}cosh(\varphi) + \boldsymbol{s}sinh(\varphi) \tag{54}$$

$$\eta_0 = -\frac{1}{\sqrt{2H}}s_0\,cosh(\varphi) - \frac{1}{\sqrt{2H}}r_0 sinh(\varphi) \tag{55}$$

$$\boldsymbol{\eta} = -\frac{1}{\sqrt{2H}}\boldsymbol{s}\,cosh(\varphi) - \frac{1}{\sqrt{2H}}\boldsymbol{r}sinh(\varphi) \tag{56}$$

$F_1$ is the algebraic part (factor) of the forward *LS* mapping and $F_2$ is the trigonometric part.

The 2x2 matrix representing $F_2$ is very much like a rotation matrix, and is easy to invert, resulting in mapping $G_2$.

*Proposition 4.* $F = F_2 \circ F_1$.

**Proof.** This is a simple substitution. ∎

*Definition 4. For positive energy, the inverse LS mapping is defined by*

i. *Mapping $G$*

$$\varphi = \xi_0 sinh(\varphi) + \frac{\eta_0}{\sqrt{-\eta^2}}cosh(\varphi) \tag{57}$$

$$\boldsymbol{q} = -\eta^2\left[\boldsymbol{\xi}\left(\frac{\eta_0}{\sqrt{-\eta^2}} + sinh(\varphi)\right) - \frac{\boldsymbol{\eta}}{\sqrt{-\eta^2}}(\xi_0 - cosh(\varphi))\right] \tag{58}$$

$$\boldsymbol{p} = -\frac{\boldsymbol{\xi}cosh(\varphi) + \frac{\boldsymbol{\eta}}{\sqrt{-\eta^2}}sinh(\varphi)}{\sqrt{-\eta^2}\left(1 - \xi_0 cosh(\varphi) - \frac{\eta_0}{\sqrt{-\eta^2}}sinh(\varphi)\right)} \tag{59}$$

ii. *Mapping $G_1$*

$$\boldsymbol{q} = -\frac{1}{2H}\left[\boldsymbol{s}(1 - r_0) + \boldsymbol{r}s_0\right] \tag{60}$$

$$\boldsymbol{p} = -\frac{\sqrt{2H}\boldsymbol{r}}{(1 - r_0)} \tag{61}$$

iii. *Mapping $G_2$*

This is the trigonometric part (factor) of the inverse *LS* mapping.

$$\varphi = \xi_0 sinh(\varphi) + \frac{\eta_0}{\sqrt{-\eta^2}}cosh(\varphi) \tag{62}$$

$$G_2 = \begin{pmatrix} cosh(\varphi) & \sqrt{2H}sinh(\varphi) \\ -sinh(\varphi) & -\sqrt{2H}cosh(\varphi) \end{pmatrix} \tag{63}$$

$$r_0 = \xi_0 cosh(\varphi) + \frac{\eta_0}{\sqrt{-\eta^2}}sinh(\varphi) \tag{64}$$

$$r = \xi cosh(\varphi) + \frac{\eta}{\sqrt{-\eta^2}}sinh(\varphi) \tag{65}$$

$$s_0 = -\xi_0 sinh(\varphi) - \frac{\eta_0}{\sqrt{-\eta^2}}cosh(\varphi) \tag{66}$$

$$\boldsymbol{s} = -\boldsymbol{\xi}sinh(\varphi) - \frac{\boldsymbol{\eta}}{\sqrt{-\eta^2}}cosh(\varphi) \tag{67}$$

The denominator of (59) $1 - \xi_0 cosh(\varphi) - \frac{\eta_0}{\sqrt{-\eta^2}}sinh(\varphi) \neq 0$, and (61) $1 - r_0 \neq 0$, because, following the identity (104), $-2Hq = 1 - \xi_0 cosh(\varphi) - \frac{\eta_0}{\sqrt{-\eta^2}}sinh(\varphi)$, and the mapping is only defined for $H > 0$, corresponding to the hyperbolic orbits of the Kepler Hamiltonian, and $q \neq 0$, meaning that collision orbits are excluded.

Now we can invert mapping $F_1$ to produce $G_1$ (details in supporting information).

Proposition 5.

i.  $G = G_1 \circ G_2$.

ii. $G_1 = F_1^{-1}$

iii. $G_2 = F_2^{-1}$

iv. $G = F^{-1}$.

**Proof.** This is almost the same as for negative energy (details in supporting information). ∎

## Identities

The following identities are direct results of the definitions and can be calculated in either direction, i.e. solely based on the forward mapping or the inverse mapping (details in supporting information).

**Proposition 6.** *The following identities hold:*

$$r^2 = r_0^2 - \sum_{i=1}^3 r_i^2 = 1 \tag{68}$$

$$r_0 s_0 - \sum_{i=1}^3 r_i s_i = 0 \tag{69}$$

$$s^2 = s_0^2 - \sum_{i=1}^3 s_i^2 = -1 \tag{70}$$

$$\xi^2 = \xi_0^2 - \sum_{i=1}^3 \xi_i^2 = 1 \tag{71}$$

$$\xi_0 \eta_0 - \sum_{i=1}^3 \xi_i \eta_i = 0 \tag{72}$$

$$\eta^2 = \eta_0^2 - \sum_{i=1}^3 \eta_i^2 = -\frac{1}{2H} \tag{73}$$

$$\varphi = \sqrt{2H} \boldsymbol{q} \cdot \boldsymbol{p} = -s_0 = \xi_0 sinh(\varphi) + \frac{\eta_0}{\sqrt{-\eta^2}} cosh(\varphi) \tag{74}$$

$$-2Hq = -p^2 q + 2 = 1 - r_0 = 1 - \xi_0 cosh(\varphi) - \frac{\eta_0}{\sqrt{-\eta^2}} sinh(\varphi) \tag{75}$$

$$L_i = q_j p_k - q_k p_j = -\frac{1}{\sqrt{2H}} \left( r_j s_k - r_k s_j \right) = \xi_j \eta_k - \xi_k \eta_j \tag{76}$$

$$M_i = \frac{1}{\sqrt{2H}} \left( \frac{q_i}{q} + p_i(\boldsymbol{q} \cdot \boldsymbol{p}) - q_i p^2 \right) = -\frac{1}{\sqrt{2H}} (r_0 s_i - r_i s_0) = \xi_0 \eta_i - \xi_i \eta_0 \tag{77}$$

$$H = \frac{p^2}{2} - \frac{1}{q} = H = -\frac{1}{2\eta^2} \tag{78}$$

$$\frac{p^2}{2} - \frac{1}{q} + q_2 p_1 - q_1 p_2 = H - \frac{1}{\sqrt{2H}} (r_2 s_1 - r_1 s_2) = -\frac{1}{2\eta^2} + \xi_2 \eta_1 - \xi_1 \eta_2 \tag{79}$$

**Proof.** These are all direct results the definitions and can be calculated in either direction, i.e. solely based on the forward mapping or the inverse mapping (details in supporting information). ∎

The identities (76)–(79) show that the angular momentum and the rescaled Runge-Lenz vector both map to angular momentum, the Kepler Hamiltonian maps to the Delaunay Hamiltonian, and the Hamiltonian of the planar rotating Kepler problem keeps it form.

## Applications

**Investigate Kepler function.**   The inverse *LS* mapping begins by solving a transcendental equation (Eq (22))

$$\varphi = \xi_0 sin(\varphi) - \frac{\eta_0}{\eta} cos(\varphi)$$

We refer to this as the "generalized Kepler equation", and define what we call the "Kepler function":

$$KF : [-1,1]^2 \to \mathbb{R} : (x, y) \mapsto solution\ of\ \varphi = x\sin(\varphi) - y\cos(\varphi) \qquad (80)$$

Since we have no solution for the Kepler function $KF$ in closed form, we explore it using a numeric solution. Here, we numerically solve the generalized Kepler equation and calculate $\varphi$ for selected values of $x$ and $y$. Fig 2 provides an overview of the Kepler function.

The numerical solution, calculated on a grid from -1 to 1 with spacing .01, demonstrates the following values:

$$min(\varphi) = (-1.2587(1, 1))$$

$$max(\varphi) = (1.2587(1, -1))$$

$$min(\nabla(\varphi)) = (-4.9081(1, 0.01)\ -100(0.99, 0))$$

$$max(\nabla(\varphi)) = (4.9081(1, -0.01)\ -0.18667(1, -1))$$

## Calculate orbits and time of flight

This application can be thought of as a warm-up exercise for using the inverse *LS* mapping. The results of the calculation have certainly been made many times before, using a slightly different technique. The only new aspect is the way that we use the Delaunay Hamiltonian and the inverse *LS* mapping to provide a complete description of orbits.

Here, we provide software that accepts values for $q$ and $p$ and calculates the Hamiltonian, the period of the orbit, the angular momentum, the Runge-Lenz vector, the values of $\xi$ and $\eta$ created by the *LS* mapping, the true, eccentric, and mean anomaly, the Kepler elements (semi-major axis, eccentricity, longitude of ascending node $\Omega$, inclination, argument of pericenter $\omega$, and time of passage), and the Delaunay variables ($\ell, g, h, \mathfrak{L}, \mathfrak{G}, \mathfrak{H}$). Then we can convert the coordinates using the *LS* mapping, calculate the orbit, including time-of-flight, by using the Delaunay Hamiltonian, and convert back to $q$ and $p$ using the inverse *LS* mapping. Finally, we show all Delaunay variables during the orbit. Fig 3 provides an overview of the output of the software for calculating orbits.

The notation and the calculations used here are described in detail in the supporting information.

## Birkhoff conjecture for the circular restricted 3-body problem

The Birkhoff conjecture was originally formulated in 1915 in an investigation of the 3-body problem[30] and can be formulated in modern terms as the existence of a global surface of section[13]. In turn, the existence of a global surface of section reveals a lot about the orbit structure and can be used to better understand the dynamics[13].

Based on recent work, it may be possible to use holomorphic curves and confirm the Birkhoff conjecture by demonstrating that the energy hypersurface is convex[13, 31]. To do that, the Hamiltonian is mapped to $\mathbb{C}^2$ by using the inverse *LS* mapping, the Levi-Civita mapping, and a stereographic projection. For the 2-body problem, the Hamiltonian takes a simple form after applying the inverse *LS* mapping, but this is not true for the 3-body problem, so we investigate it using a numerical solution.

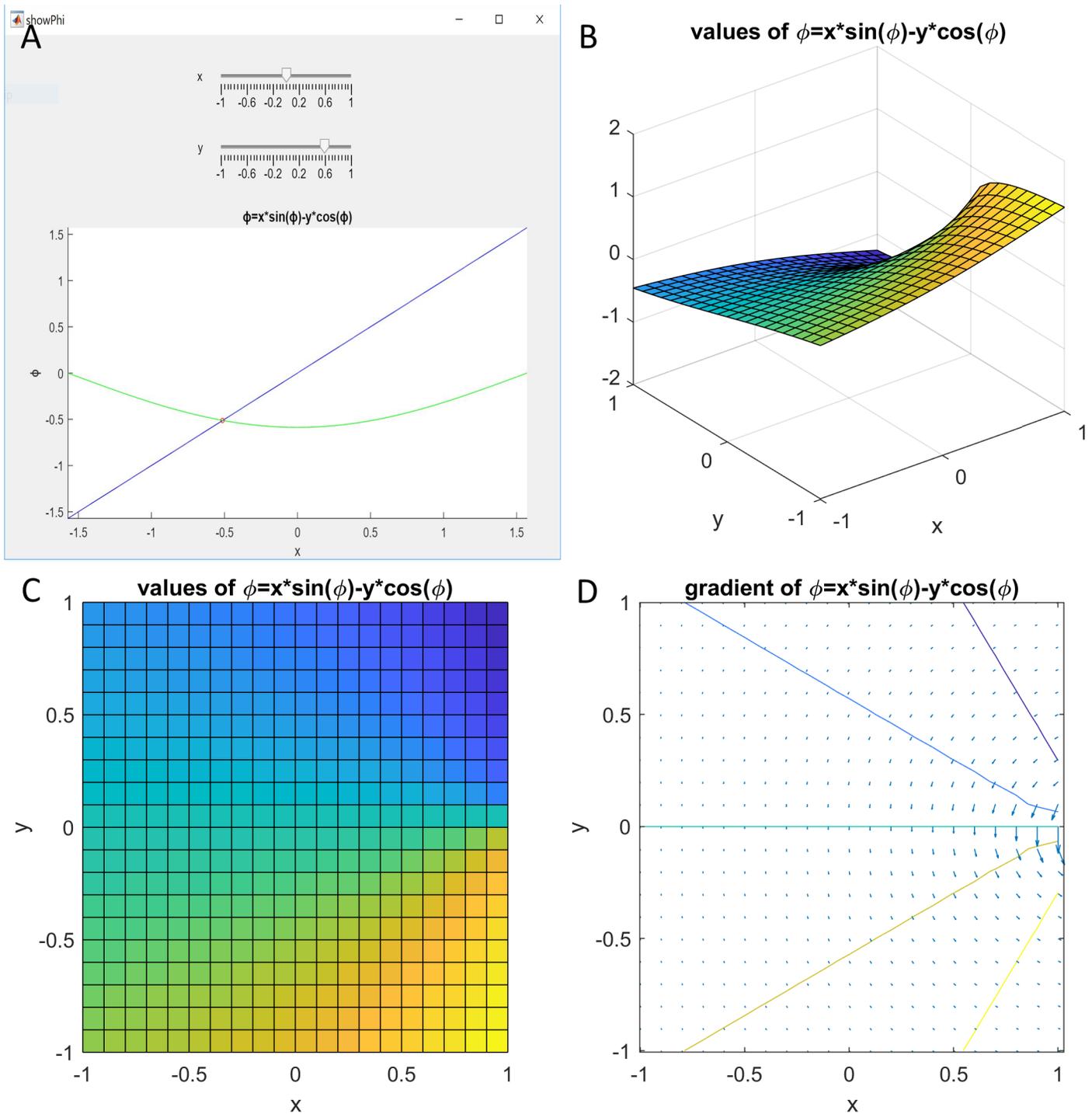

**Fig 2. Kepler function.** (A) The application allows the user to choose values for x and y, and then plots $\varphi$ (blue line), $x\sin(\varphi) - y\cos(\varphi)$ (green line), and their intersection (red circle). (B) and (C) The values of the Kepler function in two different orientations. (D) The contour and gradient.

Here we introduce the stereographic projection and the Levi-Civita mapping as used in [13]. We will use these mappings below. Fig 4 provides an overview of the stereographic projection and the Levi-Civita mapping.

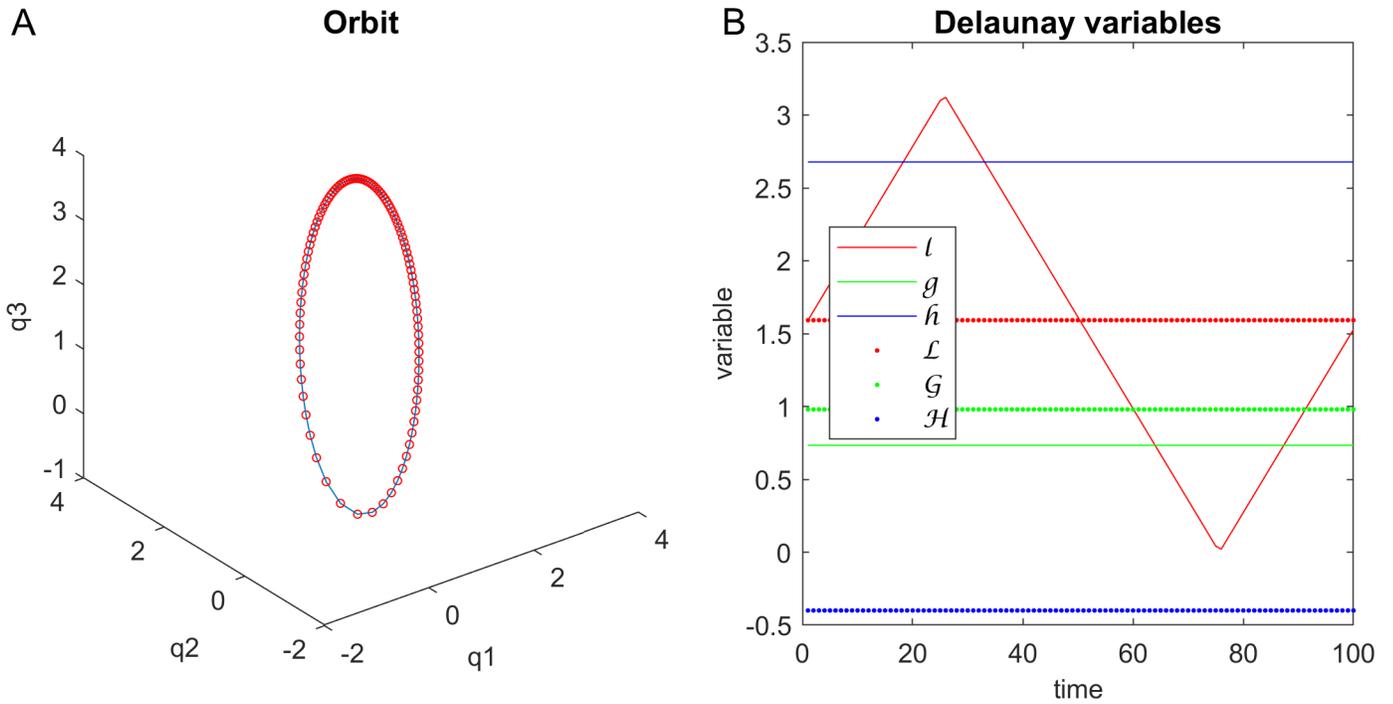

**Fig 3. Calculations and orbits.** (A) A typical orbit, with points closer where the planet is moving slower. (B) The six Delaunay variables, all of which, except for the mean eccentricity, remain constant over the orbit. The user interface is shown in S1 Fig.

https://doi.org/10.1371/journal.pone.0203821.g003

*Definition 5.* Stereographic projection.

i. Mapping SPF (stereographic projection forward)

$$\xi_0 = \frac{\|x\|^2 - 1}{\|x\|^2 + 1} \tag{81}$$

$$\xi_1 = \left(\frac{2}{\|x\|^2 + 1}\right) Re(x) = \left(\frac{2}{\|x\|^2 + 1}\right) x_1 \tag{82}$$

$$\xi_2 = \left(\frac{2}{\|x\|^2 + 1}\right) Im(x) = \left(\frac{2}{\|x\|^2 + 1}\right) x_2 \tag{83}$$

$$\eta_0 = Re(\bar{x}y) = x_1 y_1 + x_2 y_2 \tag{84}$$

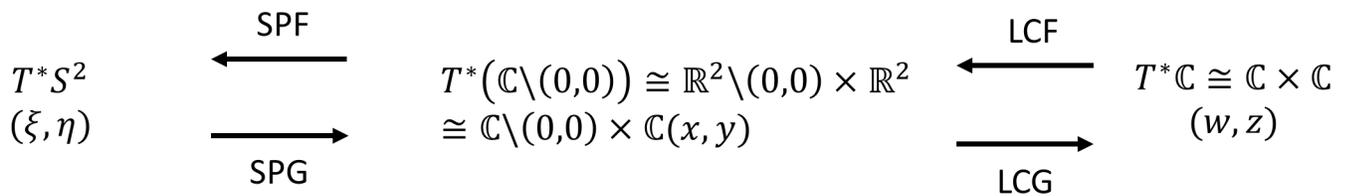

**Fig 4. Stereographic projection and Levi-Civita mapping.** The stereographic projection (SPF forward and SPG inverse) and the Levi-Civita mapping (LCF forward and LCG inverse).

https://doi.org/10.1371/journal.pone.0203821.g004

$$\eta_1 = \left(\frac{\|x\|^2 + 1}{2}\right) Re(y) - Re(\bar{x}y) Re(x) = \left(\frac{\|x\|^2 + 1}{2}\right) y_1 - (x_1 y_1 + x_2 y_2) x_1 \tag{85}$$

$$\eta_2 = \left(\frac{\|x\|^2 + 1}{2}\right) Im(y) - Re(\bar{x}y) Im(x) = \left(\frac{\|x\|^2 + 1}{2}\right) y_2 - (x_1 y_1 + x_2 y_2) x_2 \tag{86}$$

ii. Mapping SPG (stereographic projection inverse)

$$\|x\|^2 = \frac{1 + \xi_0}{1 - \xi_0} \tag{87}$$

$$x = \frac{\xi_1 + i\xi_2}{1 - \xi_0} \tag{88}$$

$$x_1 = \frac{\xi_1}{1 - \xi_0} \tag{89}$$

$$x_2 = \frac{\xi_2}{1 - \xi_0} \tag{90}$$

$$y = \eta_0(\xi_1 + i\xi_2) + (1 - \xi_0)(\eta_1 + i\eta_2) \tag{91}$$

$$y_1 = \eta_0 \xi_1 + (1 - \xi_0)\eta_1 \tag{92}$$

$$y_2 = \eta_0 \xi_2 + (1 - \xi_0)\eta_2 \tag{93}$$

The stereographic projection maps $T^*\mathbb{R}^2$ to $T^*S^2$ and the inverse projection maps to $T^*S^2$ back to $T^*\mathbb{R}^2$.

**Definition 6.** *Levi-Civita mapping.*

(ii) Mapping LCF (Levi-Civita forward)

$$x = \frac{w}{\bar{z}} \tag{94}$$

$$y = 2z^2 \tag{95}$$

(ii) Mapping LCG (Levi-Civita inverse)

$$w = \mp \frac{1}{\sqrt{2}} x \overline{\sqrt{y}} \tag{96}$$

$$z = \mp \frac{1}{\sqrt{2}} \sqrt{y} \tag{97}$$

The forward Levi-Civita mapping, in our notation, maps $T^*\mathbb{C}^2$ to $T^*\mathbb{R}^2$ and the inverse Levi-Civita mapping, maps $T^*\mathbb{R}^2$ to $T^*\mathbb{C}^2$.

## Identities

**Proposition 7.** *The following identities hold*:

$$\xi^2 = \xi_0^2 + \xi_1^2 + \xi_2^2 = 1 \tag{98}$$

$$\xi \cdot \eta = \xi_0 \eta_0 + \xi_1 \eta_1 + \xi_2 \eta_2 = 0 \tag{99}$$

$$\eta^2 = \eta_0^2 + \eta_1^2 + \eta_2^2 = \left( \frac{\|x\|^2 + 1}{2} \right)^2 \|y\|^2 = (\|w\|^2 + \|z\|^2)^2 \tag{100}$$

$$\xi_1 \eta_2 - \xi_2 \eta_1 = x_1 y_2 - x_2 y_1 = 2(w_1 z_2 - w_2 z_1) \tag{101}$$

$$H_K = -\frac{1}{2\eta^2} = \frac{2}{(\|x\|^2 + 1)^2 \|y\|^2} = -\frac{1}{2(\|w\|^2 + \|z\|^2)^2} \tag{102}$$

$$H_R = \frac{\|p\|^2}{2} - \frac{1}{\|q\|} + (p_1 q_2 - p_2 q_1)$$

$$= -\frac{1}{2\eta^2} + \xi_1 \eta_2 - \xi_2 \eta_1$$

$$= \frac{2}{(\|x\|^2 + 1)^2 \|y\|^2} + x_1 y_2 - x_2 y_1$$

$$= -\frac{1}{2(\|w\|^2 + \|z\|^2)^2} + 2(w_1 z_2 - w_2 z_1) \tag{103}$$

$$H_C = \frac{\|p\|^2}{2} - \frac{\mu}{\|q - m\|} - \frac{1 - \mu}{\|q - e\|} + (p_1 q_2 - p_2 q_1) \tag{104}$$

where $e = (-\mu, 0), m = (1 - \mu, 0)$ represent the coordinates of the earth and moon, resulting in

$$H_C = \frac{\|p\|^2}{2} - \frac{\mu}{\sqrt{(q_1 + \mu - 1)^2 + q_2{}^2}} - \frac{1 - \mu}{\sqrt{(q_1 + \mu)^2 + q_2{}^2}} + (p_1 q_2 - p_2 q_1) \tag{105}$$

Identity (101) is the angular momentum, (102) is the Kepler Hamiltonian, (103) is the Hamiltonian of the planar rotating Kepler problem [13]. In (104), for the (planar) circular restricted 3-body problem in a rotating frame, we use the Hamiltonian as defined in [30], except that the sign of the angular momentum is reversed in order to match the planar rotating Kepler problem.

**Proof.** These are all direct results of the definitions and can be calculated in either direction, i.e. solely based on the forward mapping or the inverse mapping (details in supporting information). ∎

Ultimately, we would like to examine $H_X = H_C \circ T_\mu \circ LSG \circ SPF \circ LCF$, but we observe that, when we calculate $H_C \circ T_\mu \circ LSG$, the expressions for $sin(\varphi)$ and $cos(\varphi)$ do not cancel out, leaving us with an expression that we can only solve numerically. In addition, we have added a small coordinate transformation, $T_\mu$, that shifts $q_1$ by a value of $\mu$. This is necessary because the $LS$ mapping regularizes for the collision point, $q = 0$, but in the 3-body problem, we need to regularize for the heavy primary, $q = \mu$.

Now we want to find the curvature of the tangential component of the Hessian using the method of [13]. If this curvature is always positive, we know that the energy hyperspace is convex. Specifically, the following steps are involved:

For a fixed value of the energy $c$, we first find the points of energy hyperspace $\Sigma_c = H^{-1}(c)$. This level set, or energy hyperspace, is 4-dimensional, and we consider its surface, which is a 3-dimensional hypersurface. Now we want to find the Gauss-Kronecker curvature tangential to this hypersurface. Let $G$ denote the gradient of the Hamiltonian, i.e.

$$G = \nabla H \tag{106}$$

then G is a 4-vector and we can think of it as a quaternion

$$G = (G_1, G_2 i, G_3 j, G_4 k)^T \tag{107}$$

Since G is perpendicular to the energy hypersurface, the following 3 vectors form a basis for the tangential component of G:

$$G_{tan} = (Gi \; Gj \; Gk) = \begin{pmatrix} -G_2 & -G_3 & -G_4 \\ G_1 & -G_4 & G_3 \\ G_4 & G_1 & -G_2 \\ -G_3 & G_2 & G_1 \end{pmatrix} \tag{108}$$

Now we can find the curvature:

$$C_K = det(G_{tan}{}^T * Hess(H) * G_{tan}) \tag{109}$$

For the Kepler Hamiltonian, these expressions can be calculated in closed form.

**Proposition 8.** *For the Kepler Hamiltonian,*

i. *The tangential curvature of the Hamiltonian is*

$$C_K = \frac{512}{(\|w\|^2 + \|z\|^2)^{24}}$$

ii. *The tangential curvature is always positive.*

iii. *The bounded component of the energy hypersurface $\Sigma_c = H^{-1}(c)$ is convex.*

**Proof.** (i) The Kepler Hamiltonian (102) is

$$H_K = -\frac{1}{2(\|w\|^2 + \|z\|^2)^2} \tag{110}$$

The gradient of the Hamiltonian is

$$G_K = \nabla H_K = \frac{2}{(\|w\|^2 + \|z\|^2)^3} \begin{pmatrix} w_1 \\ w_2 \\ z_1 \\ z_2 \end{pmatrix} \tag{111}$$

The Hessian of the Hamiltonian is

$$Hess(H_K) = \frac{2}{(\|w\|^2 + \|z\|^2)^3} \begin{pmatrix} 1 & 0 & 0 & 0 \\ 0 & 1 & 0 & 0 \\ 0 & 0 & 1 & 0 \\ 0 & 0 & 0 & 1 \end{pmatrix}$$
$$- \frac{12}{(\|w\|^2 + \|z\|^2)^4} \begin{pmatrix} w_1^2 & w_1 w_2 & w_1 z_1 & w_1 z_2 \\ w_1 w_2 & w_2^2 & w_2 z_1 & w_2 z_2 \\ w_1 z_1 & w_2 z_1 & z_1^2 & z_1 z_2 \\ w_1 z_2 & w_2 z_2 & z_1 z_2 & z_2^2 \end{pmatrix} \tag{112}$$

The tangential component of the gradient is

$$G_{Ktan} = \frac{2}{(\|w\|^2 + \|z\|^2)^3} \begin{pmatrix} -w_2 & -z_1 & -z_2 \\ w_1 & -z_2 & z_1 \\ z_2 & w_1 & -w_2 \\ -z_1 & w_2 & w_1 \end{pmatrix} \tag{113}$$

The curvature is

$$C_K = det(G_{Ktan}{}^T * Hess(H) * G_{Ktan}) = \frac{512}{(\|w\|^2 + \|z\|^2)^{24}} \tag{114}$$

(ii)-(iii) For the Kepler Hamiltonian, we can see that the curvature is positive for all values of the energy, so the energy hyperspace $H^{-1}(c)$ is always convex. ∎

For the planar rotating Kepler problem, these expressions can also be calculated in closed form, but the expression for the curvature is quite long. We have included all details in the supporting information. The paper on which we have modeled this calculation [13] calculates the curvature of a mapping based on the Hamiltonian and a fixed value of $c$, whereas we calculate the curvature of a mapping based directly on the Hamiltonian. The final result, which states that the energy hyperspace $H^{-1}(c)$ is always convex, is the same, but the proof of positive curvature is different, meaning that we have provided an alternate proof of the result.

**Proposition 9.** *For the planar rotating Hamiltonian, and* $c = H_R \le -1.5$,

i. *The tangential curvature of the Hamiltonian is*

$$C_{Rt} = \left(\frac{512}{X^{24}}\right) C_1 C_2 C_3 C_4 C_5 C_6$$

where

$$C_1 = X - 1$$

$$C_2 = X + 1$$

$$C_3 = X^2 + X + 1$$

$$C_4 = X^2 - X + 1$$

$$C_5 = -6L^2X^4 + LX^8 - LX^2 + 7X^6 - 1$$

$$C_6 = (X^6 + 2LX^2 + 1)^2$$

$$X = w_1^2 + w_2^2 + z_1^2 + z_2^2$$

$$L = 2(w_1z_2 - w_2z_1)$$

ii. *The tangential curvature is always positive.*

iii. *The bounded component of the energy hypersurface $\Sigma_c = H^{-1}(c)$ is convex.*

**Proof.** (i)-(ii) These are very long calculations, and the details are provided in the supporting information. (iii) is a direct result of (ii). In the proof, we make use of the assumptions $|q| \leq 1$, which holds for the bounded component of the Hill region and is equivalent to $H_K \leq -.5$, and $c = H_R \leq -1.5$. ∎

For the circular restricted 3-body problem in a rotating frame, it is not possible to calculate the curvature in closed form, because the Kepler function does not cancel out of the inverse *LS* mapping. Numerical calculations show positive values of the curvature within certain limitations.

**Observation 10.** *For the circular restricted 3-body Hamiltonian, and $c \leq$ the energy of the Lagrangian point $L_1$, numerical calculations based on the inverse LS mapping provide evidence that*

i. *The tangential curvature is positive with limitations stated in (iii).*

ii. *The bounded component of the energy hypersurface $\Sigma_c = H^{-1}(c)$ is convex with limitations stated in (iii).*

iii. *If rdHP (relative distance to the heavy primary) is the distance of the secondary object to the heavy primary divided by the distance of the Lagrangian point $L_1$ to the heavy primary, then (i) and (ii) are true in the following regions*

   a. *for $\mu = 0$ (planar rotating Kepler problem), $0 \leq rdHP < 1$.*

   b. *for $\mu = 3.277 * 10^{-7}$(Sun-Mars), $0.005 < rdHP < 1$.*

   c. *for $\mu = 3.003 * 10^{-6}$(Sun-Earth), $0.01 < rdHP < 1$.*

   d. *for $\mu = 9.536 * 10^{-4}$(Sun-Jupiter), $0.11 < rdHP < 1$.*

   e. *for $\mu = 1.216 * 10^{-2}$(Earth-Moon), $0.15 < rdHP < 1$.*

f. *for μ = 0.1, 0.3 < rdHP < 1 and c ≤ −1.8.*

g. *for μ = 0.5, 0.3 < rdHP < 1 and c ≤ −2.0.*

**Details.** The details consist of two sets of numerical calculations of the tangential curvature. In one method, we create a grid of points in $\mathbb{C}^2$ and calculate the curvature for each one. The validity of the method is limited by the fact that we can never know when we have calculated enough points. We also record some other values (e.g. $q$) in hopes of seeing a pattern.

The other method consists of calculating a constrained minimum of the curvature, where the constraints are "distance to the heavy primary is less than or equal to the distance of the heavy primary to $L_1$", and "$c$ is less than or equal to the energy of $L_1$". This provides a numerical proof that the curvature in the basin of attraction of the local minimum is always positive. The validity of this method is limited by the fact that we can never know if we have found all local minima, i.e. a global minimum.

The details of the calculations, showing some patterns that arise, are included in the supporting information.

The loss of positive curvature for points close to the heavy primary is discussed in more detail in the supporting information. It appears to arise from the fact that points close to the heavy primary can have a wide range of energy values, leading to a higher probability of a larger curvature. For increasing values of, this happens for increasing values of $q$, so it is not simply an issue with getting too close to the singular point, where the software used to calculate gradient and Hessian might become erratic.

**Observation 11.** *For the circular restricted 3-body Hamiltonian, and c ≤ the energy of the Lagrangian point $L_1$, numerical calculations based on the inverse LS mapping provide evidence that the Birkhoff conjecture is true.*

**Details.** We follow the argumentation of Frauenfelder et al [13, 31]. In theorem 1.3, Hofer, Wysocki, and Zehnder [32] used holomorphic curves to show the existence of a global surface of section under the condition that the energy hypersurface is contact and dynamically convex, i.e. the Conley-Zehnder indices of all periodic orbits are greater than or equal to three. Moreover, they showed that if the energy hypersurface admits a convex embedding it is dynamically convex. Later, in Theorem A, Albers, Frauenfelder, van Koert, and Paternain [33] proved the contact condition, so that only dynamical convexity remains to be checked. More details can be found in [31]. As discussed in [13], numerical experiments of Otto van Koert, based on the Levi-Civita embedding, show that nonconvex points arise close to the Lagrange points. In contrast with that, numerical evidence for Proposition 10 (above) shows that, for Sun-Mars, Sun-Earth and Sun-Jupiter, the only nonconvex points are close to collisions. That there are nonconvex points close to collisions is due to the smoothness issue of the Ligon-Schaaf map at collisions. This implies that the Conley-Zehnder indices of periodic orbits which do not come close to the sun have Conley-Zehnder index greater than or equal to three. Combining both embeddings and believing the numerical evidence then implies that the only periodic orbits which could have Conley-Zehnder index less than three are periodic orbits which come close to the sun as well as to the Lagrange point. Probably, such orbits have higher actions than the retrograde periodic orbit and therefore do not obstruct the Birkhoff conjecture, telling us that the retrograde periodic orbit bounds a global surface of section.

With this, the validity of the Birkhoff conjecture, i.e. the existence of a global surface of section, can be used to study the dynamics of orbits, as discussed in [13]. Due to the work of Franks, Handel and Le Calvez, we now have a much deeper understanding of the nature of area-preserving disk maps [34, 35]. However, specific properties of the return map of the circular restricted 3-body problem are still open and these form in fact one of the ultimate goals beyond the Birkhoff conjecture.

## Supporting information

**S1 Document. Supporting information to the symmetry of the Kepler problem inverse Ligon-Schaaf mapping and the Birkhoff conjecture.** This document contains proofs and other detailed calculations.
(PDF)

**S1 Fig. LS orbit software.** Example Application (software used for Fig 3.)
(TIFF)

**S2 Fig. $C_5$ constraints and plots.** Part of proof of proposition 9.
(TIFF)

**S3 Fig. $C_5$ plots at specific points.** Part of proof of proposition 9.
(TIFF)

**S4 Fig. Sun-Mars and Sun-Earth energy and curvature.** Part of discussion of proposition 10.
(TIFF)

**S5 Fig. Sun-Jupiter energy and curvature.** Part of discussion of proposition 10.
(TIFF)

**S6 Fig. Sun-Jupiter energy and curvature 2.** Part of discussion of proposition 10.
(TIFF)

**S7 Fig. Earth-Moon energy and curvature.** Part of discussion of proposition 10.
(TIFF)

**S8 Fig. Hill regions.** Hill regions for $\mu = 0.1$ and $\mu = 0.5$. Part of discussion of proposition 10.
(TIFF)

**S9 Fig. $\mu = 0.1$ and $\mu = 0.5$ energy and curvature.** Part of discussion of proposition 10.
(TIFF)

**S1 Software. MATLAB software.** This .zip file contains all of the MATLAB software used to create data and diagrams in the paper (details in supporting information).
(ZIP)

**S1 Data.** This .zip file contains data created by the software in S1 Software.
(ZIP)

## Acknowledgments

I would like to thank Urs Frauenfelder for many useful discussions without which this project would not have been possible. My thanks also go to Lei Zhao for contacting me via social media and suggesting a meeting that led to this project.

## Author Contributions

**Conceptualization:** Thomas Sumner Ligon.

**Investigation:** Thomas Sumner Ligon.

**Methodology:** Thomas Sumner Ligon.

**Software:** Thomas Sumner Ligon.

**Visualization:** Thomas Sumner Ligon.

**Writing – original draft:** Thomas Sumner Ligon.

**Writing – review & editing:** Thomas Sumner Ligon.